\documentclass[12pt,psamsfonts]{article}
\usepackage{amsmath}
\usepackage{amssymb}
\usepackage{amsthm}

\def\ui{{\underline i}}

\def\ua{{\underline a}}

\def\depth{\operatorname{depth}}

\def\gin{\operatorname{Gin}}

\newcommand{\m}{\mathfrak m}
\newcommand{\n}{\mathfrak n}
\newcommand{\M}{\mathfrak m}

\def\gin{\mathop{\kern0pt\fam0Gin}\nolimits}
\def\gr{\mathop{\kern0pt\fam0gr}\nolimits}
\def\Tor{\mathop{\kern0pt\fam0Tor}\nolimits}
\def\depth{\mathop{\kern0pt\fam0depth}\nolimits}
\def\ini{\mathop{\kern0pt\fam0in}\nolimits}
\def\gcd{\mathop{\kern0pt\fam0gcd}\nolimits}
\def\Grass{\mathop{\kern0pt\fam0Grass}\nolimits}

\def\gin{\operatorname{gin}}

\newtheorem{theorem}{Theorem}[section]

\newtheorem{corollary}[theorem]{Corollary}
\newtheorem{proposition}[theorem]{Proposition}
\newtheorem{remark0}[theorem]{Remark}
\newtheorem{example0}[theorem]{Example}

\newenvironment{example}{\begin{example0}\rm}{\end{example0}}
\newenvironment{remark}{\begin{remark0}\rm}{\end{remark0}}

\newcommand{\propref}[1]{Proposition~\ref{#1}}
\newcommand{\thmref}[1]{Theorem~\ref{#1}}

\newcommand{\corref}[1]{Corollary~\ref{#1}}

\title{  \bf   Isomorphism classes of short Gorenstein  local rings via Macaulay's  inverse system
}
\author{\large   J. Elias
\thanks{Partially supported by  MTM2007-67493, Acci\'{o}n Integrada Espa\~{n}a-Italia 07-09}
\and \large M. E. Rossi
\thanks{Partially supported by M.I.U.R.: PRIN 07-09, Azione Integrata Italia-Spagna 07-09 }
}

\date{  }

\begin{document}

\maketitle

\medskip
\begin{abstract}
In this paper we study the isomorphism classes of    Artinian Gorenstein local rings  with socle degree three  by means of Macaulay's inverse system.
We prove that their  classification  is equivalent to the  projective classification of the
hypersurfaces of $\mathbb P ^{n }$ of degree three. This is an unexpected result because it  reduces the study of this class of local rings to the homogeneous case.    The result has applications   in problems concerning the  punctual Hilbert scheme $ Hilb_d (\mathbb P^n) $
   and in relation to the problem of the rationality of the Poincar\'e series of local rings. 
\end{abstract}

\section{Introduction}

The classification,  up  analytic isomorphism,    of  Artinian local algebras    plays an important role  in commutative algebra and algebraic geometry.    Among other examples, the study of the irreducibility and   the
smoothness of the punctual Hilbert scheme $ Hilb_d (\mathbb P^n) $
parameterizing zero-dimensional subschemes of  fixed degree $d $ in $\mathbb P^n,    $ is strictly related to the structure of the Artinian local algebras of multiplicity $d. $ In particular,  because the locus of points of  $ Hilb_d (\mathbb P^n) $ corresponding to Gorenstein subschemes is   open and non-empty, often  one can restrict  the study to schemes which are the spectrum of   Artinian Gorenstein local algebras,  see  for example \cite{CEVV08},  \cite{CN07},  \cite{CN09}, \cite{EV}, \cite{I}, \cite{Po}. \par \noindent  In this paper we present  structure theorems for Artinian Gorenstein local rings $(A, \m)$ such that $\m^4=0.  $     Over the last   years the  class under investigation has  emerged as a testing ground for properties of infinite free resolutions (see \cite{AIS}).  Rings of low socle degree have  been studied indeed in relation to the problem of the rationality of the Poincar\'e series. We recall that Anick in \cite{A}  constructed the first example of a local ring  for which the Poincar\'e series is
 irrational.  Notice that  it satisfies $\m^3=0.$ Sj\"odin proved that all Artinian Gorenstein   rings with $\m^3 =0$  have a rational Poincar\'e series,
\cite{Soj79}, B{\o}gvad proved that there exist    Artinian Gorenstein  local rings with $\m^4=0 $   and irrational Poincar\'e series, \cite{Bog83}. However  Henriques and \c{S}ega recently proved that,  assuming the existence of exact zero-divisors,  the Poinca\'re series of an Artinian  Gorenstein local ring of socle degree three is rational,  \cite{HS09}.
Our approach   will   carry   the problem over  graded  algebras where  Henriques and \c{S}ega's result  was already considered in  \cite{CRV01} (see   Corollary  \ref{ines}). 
\vskip 2mm

In the main result of this paper we prove that the classification of   Artinian Gorenstein  local  algebras  $(A, \m) $ with
$\m^4=0 $  is equivalent  to the projective
classification of the  hypersurfaces of $\mathbb P^{n }$ of  degree three (see Theorem \ref{MAIN} and Corollary \ref{isomnp}).  The key point of the   paper is to prove  that  any     Artinian Gorenstein  local ring  with
Hilbert function $\{1, n, n, 1\}$ is canonically isomorphic to its associated graded ring  (see  \thmref{n}), hence it is canonically graded.  This is an unexpected result and   it  is false in general,  even if the Hilbert function is symmetric, for instance  if the Hilbert function is $\{1, 2, 2, 2, 1\}. $   The results of Section 3. concerning the balanced Hilbert function $\{1, n, n, 1\} $ will be extended in Section 4.  to the general situation of an Artinian Gorenstein local ring  with  socle degree three.
 As a consequence we  will take advantage of  the literature concerning graded  Artinian Gorenstein rings with socle degree $3$ for getting results in the local case.
Important classes of projective  varieties have   Artinian reduction  of low socle degree. With regard to this, we recall   that the homogeneous coordinate ring of the canonical curve of genus $g$ is a Gorenstein domain of dimension two with $h$-vector $\{1, g-2, g-2, 1\}, $  hence its Artinian reduction is  still Gorenstein with socle degree $3. $   \vskip 2mm

From  classical  results on the projective classification of homogeneous cubics in $\mathbb P^{r}$ with $r \le 3, $ we may deduce the classification of Artinian Gorenstein local rings with Hilbert function $\{1, m, n, 1\} $ where $n \le 4 $ and a nice geometric description can be given,   see \corref{iso} and \corref{V(F)}.     This  classification  allows us to find  geometric interpretations  of the singularity of  $ Hilb_d (\mathbb P^n) $ for low degrees $d.$  In particular we extend  recent   results by Casnati and Notari (see  \cite{CN09}, Theorem 4.1)  and by Cartwright et al.  in \cite{CEVV08}. 
  \par \noindent

 \vskip 2mm

A central tool in this paper is  Macaulay's inverse system (see  \cite{Mac16}) which
    establishes  a one-to-one correspondence between Artinian Gorenstein  algebras
and suitable  polynomials.
This classical correspondence has been deeply studied  in the homogeneous case, among other authors, by Iarrobino  in a long series of papers, see for example:  \cite{Iar84},  \cite{Iar94}, \cite{IK99}.
Notice that from a categorical point of view, Macaulay's correspondence is a particular case of Matlis duality, see \cite{BS98},  \cite{Nor72b} and
\cite {SV72},  Theorem 5.2. \par \noindent
  In  the local case the situation is more complicate than in the graded one, nevertheless
the classification    under the action of isomorphisms
of algebras can be  translated to the classification of the inverse system polynomials under the action of a group of transformations explicitly described. Emsalem devoted  \cite{Ems78}, Section C,   to present an  approach to the classification of
  Artinian Gorenstein local algebras by means of their inverse system. 
In Section 2.  we collect   results  spread over different papers (see  \cite{Ems78}, \cite{Iar84},  \cite{Iar94}, etc.)   by presenting  also explicit  methods.
The  strategy of classifying inverse system polynomials instead of   the corresponding   Artinian Gorenstein local  algebras  has at least two advantages: first,  we deal with one  polynomial instead of a   system  of
generators of the defining ideals  (as in   \cite{CN07}, \cite{CN09},  \cite{EV})   and,  second,  we may  perform   effective computations, often reduced to a linear  algebra problem  (see  \propref{EmsalemGor}).  We hope that this approach will be useful in studying numerical invariants of local Gorenstein singularities.

\bigskip

\noindent {\sc Acknowledgments.} The authors are grateful to A. Conca and   G. Valla for   useful comments and remarks regarding this work.

\bigskip
\section{Inverse System for Artinian local rings}

\bigskip
\noindent

Through the paper we are assuming that the basic field $K$ is algebraically closed
of characteristic zero.
Let $R=K[\![x_1,\dots x_n]\!]$ be the ring of the formal series with maximal ideal $\mathcal M =(x_1,\cdots,x_n)$
and let  $P=K[y_1,\dots,y_n]  $ be  a  polynomial ring.
Let $s $ be a positive integer, we denote by $ P_{\le s} $ the set of polynomials of degree $ \le s.$

Let  $ I$ be an ideal of $R$ such that $A= R/I $ has finite length.
Let $\m= {\mathcal M} +I $ be  the maximal ideal of $A,  $       then the socle of $A $ is the colon ideal $ 0:   \m.$
The socle-degree $ s $ of $ A $ is the largest integer for which $  \m^s \neq 0. $
$ A $ is Gorenstein if and only if  dim$_K (0:  \m ) =1.$ From now the local ring $(A, \m) $ will be  a quotient of  $R.$

Consider   $gr_{\M}(A) := \bigoplus_{i \ge 0} \M^i/\M^{i+1} $   the associated graded ring of $A=R/I.$
It is well known that
 $gr_\M(A) \simeq K[x_1,\dots x_n]/I^* $ where $I^* $ is the  homogeneous ideal of $K[x_1,\dots x_n]$ generated by the initial forms of the elements of $I.$
The Hilbert function of $A$ is by definition the Hilbert function of $gr_\M(A)$, i.e.
$$
HF_A(i) =   dim_K \left(\frac{\M^i}{\M^{i+1}}\right).
$$
\bigskip

In this section we collect the main facts and notations concerning  Macaulay's inverse system in the  study and classification  of  Artinian  local rings $(A, \m).$
The reader should refer to \cite{Ems78} and \cite{Iar94} for an extended treatment.   The graded case is much more understood   than the local case and several interesting papers have been written (see for example Chapter IV of \cite{Mac16} and in \cite{IK99}).

\vskip 3mm It is known that
$P$ has a structure of   $R$-module by means the following action
$$
\begin{array}{ cccc}
\circ: & R  \times P  &\longrightarrow &  P   \\
                       &       (f , g) & \to  &  f  \circ g = f (\partial_{y_1}, \dots, \partial_{y_n})(g)
\end{array}
$$
where $  \partial_{y_i} $ denotes the partial derivative with respect to $y_i.$
If  we  denote by $ x^{\alpha}= x_1^{\alpha_1} \cdots x_n^{\alpha_n} $ and
$ y^{\beta}= y_1^{\beta_1} \cdots y_n^{\beta_n} $ then
$$
x^{\alpha} \circ y^{\beta} =
\left\{
\begin{array}{ll}
\binom{\beta}{\alpha}\; y^{ \beta-  \alpha}  & \text{ if }  \beta_i \ge \alpha_i \text{ for } i=1,\cdots, n\\ \\
0 & \text{ otherwise}
\end{array}
\right.
$$
where $\binom{\beta}{\alpha}=\prod_{i=1}^n  \binom{\beta_i}{\alpha_i}$.
The action for polynomials is defined    bilinearly  from that for monomials.
We remark that  for every $f, h \in R $ and $g \in P, $ then $(f h) \circ g = f \circ ( h \circ g)$, and that
we have  $ \mathcal M^{s+1} \circ g =0 $ if and only if $ g \in P_{\le s} $.

\noindent Let $S \subseteq P $ be a set of polynomials, in the following we will  denote by $\langle S\rangle_K$  the $K$-vector space generated  by $S$,  and by $\langle S\rangle_R $ the $R$-submodule of $P$ generated by the set $S$, i.e. the $K$-vector space generated  by $S$ and   the corresponding derivatives of all  orders. Starting from $\circ $ we consider the exact pairing of $K$-vector spaces:
$$
\begin{array}{ cccc}
\langle\ , \ \rangle : & R  \times P  &\longrightarrow &   K   \\
                       &       (f , g) & \to  & ( f  \circ g ) (0)
\end{array}
$$Recall that it  gives an  explicit description  of the injective envelope of the $R$-module  $K$
 in terms of  the polynomial ring  $ P$ with the structure of  $R$-module via the action
 $\circ$ (see \cite{Ems78} Section B, Lemme and Proposition 1;  \cite{Iar94} page 9,
 \cite{Nor72b}). For any ideal  $I\subset R$ we define the following $ R$-submodule of $P: $
$$
 {I^{\perp}}:=\{g\in P\ |\  \langle f, g \rangle = 0 \ \ \forall f  \in I \ \}= \{ g \in P\ |\  I \circ g = 0    \}.
$$
If $I = {\mathcal M}^{s+1}$ then $ I^{\perp} $ coincides with $P_{\le s}.$
 In  general if   $R/I$ has socle-degree  $s$, then $  {I^{\perp}} $ is generated by polynomials  of degree $\le s. $
Conversely, for every $R$-submodule $M$ of $P$ we define
$$
Ann_R(M):=\{g \in R \ |\  \langle g, f \rangle = 0 \ \ \forall f \in M  \} =
Ann_R(M) =\{g \in R \ |\  g \circ M = 0  \}
$$  which  is an ideal of $R. $
If $M$ is cyclic, that is $M  =\langle f \rangle=  R \circ f $ with $f \in P, $ then we will write $ Ann_R(f). $
Recall  that  $M  =\langle  f \rangle_R $ is the  $K$-vector space generated by the polynomial $f $ and all its derivatives of every order.
Moreover $\langle  f \rangle_R  = \langle  h \rangle_R$ if and only if there exists an invertible element $u $ of $R$ such that $ f= u \circ h. $

 Emsalem in \cite{Ems78}, Section B, Proposition 2, Iarrobino in \cite{Iar84} Theorem 4.6 B  and   \cite{Iar94} Lemma 1.2,
  proved that there exists a one-to-one correspondence between  ideals $I\subseteq R  $ such that $R/I$ is an Artinian local ring and $R$-submodules  $M$ of $P$ which are finitely generated.

\noindent Since $R/I$ and $I^{\perp} $ are   finitely generated $K$-vector spaces,  it is easy to see that the action
$\langle \ ,\ \rangle$ induces the following
isomorphism of $K$-vector spaces (see \cite{Ems78} Proposition 2 (a)): \begin{equation} \label{pairing}  (R/I)^* \simeq   {I^{\perp}}. \end{equation}

\noindent
Hence $dim_K R/I  (= \text{\ multiplicity\ of } R/I) = dim_K I^{\perp}. $ As in the graded case, it is possible to compute the Hilbert function of $A=R/I $ via the inverse system.
We define the following $K$-vector space:
\begin{equation} \label{H1}
 (I^{\perp})_i := {\frac{ I^{\perp} \cap P_{\le i} +  P_{< i}}{ P_{< i}}}
\end{equation}
then by (\ref{pairing}) we can prove
\begin{equation} \label{H2}
HF_{R/I}(i) = dim_K    (I^{\perp})_i .
\end{equation}

\bigskip
We are interested in studying the automorphisms of   Artinian local rings $A= R/I $ of socle $s.$
Given a $K$-algebra $C$ we will denote by $Aut (C) $ the group of the automorphisms as $K$-algebra and by $Aut_K(C) $ as $K$-vector space.
The  automorphisms   of $ R=K[[x_1, \dots, x_n]] $ as  $K$-algebra are well known. They act as replacement of $x_i$ by
$z_i $, $i=1,\cdots, n$, such that ${\mathcal M }= (x_1,
\dots, x_n)= (z_1, \dots, z_n). $ Actually, since $ {\mathcal M }^{s+1}  \subseteq I, $ we are interested in the   automorphisms of $R/ {\mathcal M}^{s+1} $   of $K$-algebras  induced by
the projection $\pi:R\longrightarrow R/{\mathcal M }^{s+1}.$    Clearly    Aut$(R/ {\mathcal M}^{s+1})  \subseteq $   Aut$_K (R/ {\mathcal M}^{s+1}).   $

 Let $E=\{ e_i \} $ be the canonical basis     of  $ R/\mathcal M^{s+1}  $ as $K$-vector space consisting of the standard monomials  $x^{\alpha}$ ordered by increasing degree and lexicographic order,  then the dual basis of $E$ with respect the
perfect paring $\langle\, , \rangle$ is the basis $E^*=\{ e_i^* \} $ of $P_{\le s} $  where
$$
(x^{\alpha})^* = \frac 1 {\alpha !} y^{\alpha},
$$
in fact $e_i^* (e_j)= \langle e_j ,  e_i^* \rangle=\delta_{ij}$,
where $\delta_{ij}=0$ if $i\neq j$ and $\delta_{ii}=1$.
Hence for any
$\varphi  \in  Aut_{K}(R/{\mathcal M }^{s+1}) $ we may associate   a    matrix $ M(\varphi)$   with respect the basis $E$ of size $ r= dim_K (R/{\mathcal M }^{s+1}) =   {{n+s} \choose s}.  $
Summing up we have  the following  natural sequence of  morphisms of groups: \begin{equation} \label{groups} Aut (R) \stackrel{\pi}{\longrightarrow} Aut (R/\mathcal M^{s+1} ) \stackrel{\sigma}{\longrightarrow} Aut_{K}(R/\mathcal M^{s+1} ) \stackrel{\rho_E}{\longrightarrow} Gl_r(K). \end{equation}

\noindent Given  $I $ and $ J $   ideals of $R$ such that $\mathcal M^{s+1}\subset I, J, $ there
  exists an  isomorphism of  $K$-algebras
$$
\varphi : R/I \to R/J
$$
if and only if $ \varphi  $ comes from a $K$-algebra automorphism   of $ R/{\mathcal M }^{s+1} $  sending  $I/{\mathcal M }^{s+1}$ to $J/{\mathcal M }^{s+1}.$
Dualizing
$$
\varphi^* : (R/J)^* \to (R/I)^*
$$
is an isomorphism of the  $K$-vector subspaces    $ (R/I)^* =  {I^{\perp}} $ and $ (R/J)^* =  {J^{\perp}} $
of  $P_{\le s}$.
  Hence $^t M (\varphi) $ is the matrix associated to $\varphi^* $
with respect the basis $E^*  $ of $P_{\le s}.$

\noindent We can complete (\ref{groups}) by  the following commutative diagram which helps to visualize our setting:

$$
\begin{array}{ccc}
Aut_{K}(R/\mathcal M^{s+1} ) &  \stackrel{\rho_E}{\longrightarrow} & Gl_r(K) \\
\downarrow *&  & \; \downarrow \; ^t() \\
Aut_{K}(P_{\le s} )          &  \stackrel{\rho_{E^*}}{\longrightarrow} & Gl_r(K) \\
\end{array}
$$

\noindent In the following we will be interested in the   subgroup ${\mathcal R}$ of  $Aut_{K}(P_{\le s})$    (automorphisms of $ P_{\le s} $ as $K$-vector space)
   represented     by
 the matrices $^t M (\varphi) $ of  $Gl_r(K)$ with $\varphi\in Aut (R/{\mathcal M }^{s+1}) $.

\bigskip
 \begin{theorem}\label{Emsalem}{\rm{ (\cite{Ems78} Proposition 15)}}
 The classification,  up isomorphism,   of the Artinian local algebras of multiplicity $d, $ socle degree $s$ and embedding dimension $n$ is equivalent to the classification,
 up the action of ${\mathcal R},  $ of the $K$-vector subspaces of $P_{\le s}$  of dimension $d, $ stable by derivations and containing $P_{\le 1.}$
 \end{theorem}

 \bigskip

More precise results  can be stated for Artinian Gorenstein local rings.  A  local ring  $A=R/I $ is an Artinian   Gorenstein  local ring of socle degree $s $ if and only if its dual module ${I^{\perp}}  $ is a cyclic $R$-submodule of $P$ generated by a polynomial  $F \in P $ of degree $s $ (see also \cite{Kap70}, Theorem 220,  \cite{Iar94}, Lemma 1.2.).   We will denote by $A_F $ the Gorenstein Artin algebra associated to $F \in P, $ i.e. $$A_F= R/Ann_R(F). $$

\noindent Hence each  Artinian Gorenstein local ring of socle $s$ will be equipped with a polynomial $F \in P $ of degree $s.$
The polynomial $F$ is not unique, but is determined up an unit $u$ of $R$ ($F$ can be replaced by $u \circ F$).

\vskip 3mm
Our goal is to translate the classification of the Artinian Gorenstein local rings $A=R/I $  of socle $s$ in terms of the corresponding polynomials  of degree $s $ in $P.$

\noindent Let $\varphi  \in Aut (R/ \mathcal M ^{s+1}), $  from the previous facts we have
 \begin{equation}
\varphi(A_F)=A_G \ \ \text{ if\  and \ only\  if } \  (\varphi^*)^{-1}(\langle F \rangle_R)=\langle G \rangle_R.
 \end{equation}

\noindent
It is easy to verify  that
\begin{equation} \label{matrix}
\varphi(A_F)=A_G \ \ \text{ if\  and \ only\  if } \  (\varphi^*)^{-1}(F)= u \circ  G     \ \ \text{ where} \ \   u \ \ \text{ is an unit in }  R .
\end{equation}

Let  $u $ an invertible element of $R/{\mathcal M }^{s+1}, $ the   corresponding action of
 $u$ in $P_{\le s}$ is a $K$-vector space isomorphism, so we can consider the associated
matrix $N(u)\in Gl_r(K) $ with respect the basis $E^*$.
If    $F  =  b_1 e_1^*+ \dots b_r e_r^*  \in P_{\le s}, $ then we will denote the {\it row vector} of the coefficients of the polynomial w.r.t. the basis  $E^*$ by
    $$  [F]_{E^*} = (b_1, \dots, b_r).$$

\noindent
Hence from (\ref{matrix}) we deduce the following key result:

\begin{proposition}
\label{EmsalemGor}
The  Artinian  Gorenstein local rings  $A_F$ and  $A_G$ of socle degree $s$ are isomorphic
 if and only if there is $\varphi \in Aut (R/{\mathcal M }^{s+1})$ and an invertible element $u\in R/{\mathcal M }^{s+1}$ such that
 \begin{equation} \label{vector}
[G]_{E^*}  (^t N(u)  \ M(\varphi))= [F]_{E^*}.
\end{equation}
\end{proposition}

\noindent The above result enable us to study the isomorphism classes of Artin Gorenstein algebras
in a effective computational framework.
 The strategy of this paper is to classify the Artinian  local algebras by classifying their inverse systems by means  of \eqref{vector}.

\bigskip
\noindent We remark that  Proposition \ref{EmsalemGor} extends  the well known result for which in the homogeneous case    the   isomorphisms  of   Artinian Gorenstein graded standard $K$-algebras act on the homogeneous polynomials via  the action of the linear group  $GL_n(K) $,
see \cite{IK99}, \cite{Ems78}, Section C, Proposition 17.

We say that    $F \in P=K[y_1, \dots, y_n] $ is  {\it{ non  degenerate}}  if
the embedding dimension of the Gorenstein algebra  $A_F $ is  $n$, i.e. $HF_{A_F}(1)=n$.
Hence a polynomial  $F $ of degree $s$  is non degenerate  if and only if
 the dimension of the  $K$-vector space of the derivatives of order $s-1$ of $F $ is $n$ or equivalently $Ann_R(F) $ does not contain linear elements.

In order to classify the Artinian Gorenstein local rings of given multiplicity, we need information on the admissible Hilbert functions. In the graded case,   the Hilbert function of an Artinian Gorenstein algebra is symmetric.  Little is  known about the Hilbert function in the local case. The problem comes from the fact that,  in general,   the associated algebra $G=gr_{\m}(A)$   is no longer Gorenstein.

Nevertheless  Iarrobino in \cite{Iar94} proved interesting results  concerning $G.$ Consider  a filtration  of $G $  by a descending sequence of ideals:
$$ G = C(0) \supseteq C(1) \supseteq \dots \supseteq  C(s-2)  $$
whose successive quotient
$$ Q(a)= C(a)/C(a+1) $$
are  reflexive graded $G$-modules of  socle degree $ \frac{s-a}{2}  $ ($s =$ socle-degree of $A$), see \cite{Iar94}, Theorem 1.5.  Hence $ Q(a) $ has symmetric Hilbert function. The reflexivity of $Q(a) $ as $K$-vector space is induced from the nonsingular pairing  on $A $ where the ideals $(0: \m^i) $ of $A$ correspond to $(\m^i)^{\perp}. $  In the graded case $0:\m^i =\m^{s+1-i} $ from  which follows that the  Hilbert function of $A$ is symmetric. When $A$ is not graded,   $(0: \m^i) \neq  \m^{s+1-i} $ but  the duality still gives some information.

\noindent By the reflexivity,   the $G$-modules $Q(a) $   have symmetric  Hilbert function  about $ \frac{s-a}{2}.   $  Thus the Hilbert function of $A, $ while not usually itself symmetric, is the sum of symmetric functions $ HF_a:= HF_{Q(a)}  $ with offset centers (see    Proposition 1.9 \cite{Iar94}). Notice that in general $Q(a) $ are not standard $K$-algebras, except either $Q(0) $ or the case of embedding dimension of $A$ equal to two  (\cite{Iar94}, Section 2.).
\vskip 2mm
 In particular $Q(0) = G/C(1) $ is the unique (up isomorphism) graded Gorenstein quotient of $G$ with the same socle degree $s$.  Iarrobino proved that if $HF_A(n) $ is symmetric, then $G= Q(0) $ and it is Gorenstein. The same result has also been proved in a different manner by J. Watanabe in \cite{W}. Hence (see \cite{Iar94} Proposition 1.7 and \cite{Ems78} Proposition 7)  if $A$ is Gorenstein

\medskip
{\centerline{ G is Gorenstein $\iff $ $HF_A(n) $ is symmetric $\iff $  $G=Q(0) $} }

\medskip

The $G$-module $Q(0) $ plays a crucial role and it can be computed in terms of the corresponding polynomial in the inverse system. Let $F \in P $ be a polynomial of degree  $ s $ and denote by $F_s$ the form of highest degree in $F, $ that is  $F =F_s + \dots $ terms of lower degree,  then
$$ Q(0) \simeq R/Ann_R(F_s), $$
 see  \cite{Ems78} Proposition 7 and \cite{Iar94} Lemma 1.10.

\bigskip
\section{Gorenstein Artin algebras with $HF_A=\{1,n,n,1\}$}
Let $A=R/I $ be an  Artinian Gorenstein  local ring  with Hilbert function $HF_A=\{1,n,n,1\}. $ This means that the socle degree is $3$ and $HF_A(1)=HF_A(2) =n, HF_A(0)=HF_A(3) =1.$  As we have seen, since $HF_A $ is symmetric,  then  $G= gr_{\m}(A) $ is Gorenstein.  We recall that if $R=K[[x_1, \dots,x_n]], $ then $G=K[x_1, \dots, x_n]/I^* $ where $I^*$ is the ideal generated by the initial forms of the elements of $I.$ Since $A$ is Artinian, there is the natural isomorphism between  $G$ and $R/I^*R. $   In this  section we prove that there exists an    isomorphism of local rings  between $A$ and $G $  (actually between $A$ and $R/I^*R$), see Theorem \ref{n}.

It is very rare that a local ring is isomorphic to its associated graded ring.
If $A \simeq G$ (as before specified), accordingly with  Emsalem in \cite{Ems78},  we say  that $A$ is {\it{ canonically graded}}. It is a surplus  to recall that Gorenstein local rings with symmetric Hilbert function are not necessarily canonically graded.
 The following example   shows  that we cannot extend the main result of this section  to higher
socle degrees.

 \begin{example}
 Let $A$ be an Artinian   Gorenstein local ring  with Hilbert function $HF_A=\{1,2,2,2,1\}$. Then   $A$ is isomorphic to one and only one of the following  rings:
\begin{enumerate}
\item[(a)]
$R/I$ with $I=(x^4, y^2)$ and $I^{\perp}=\langle x^3 y\rangle$, i.e. $A$ is isomorphic to its associated graded ring,
\item[(b)]
$R/I$ with $I=(x^4, -x^3 + y^2)$ and $I^{\perp}=\langle x^3 y+ y^3\rangle$.
The associated graded ring is of type $(a)$ and it is  not  isomorphic to $R/I$,
\item[(c)]
$R/I$ with $I=(x^2+y^2, y^4)$ and $I^{\perp}=\langle x y(x^2-y^2) \rangle$, i.e. $A$ is isomorphic to its associated graded ring.
\end{enumerate}
The computation can be performed by using Proposition \ref{EmsalemGor}, a different approach can be found in \cite{EV}.
 \end{example}

\bigskip
Since  $G$ is a graded Gorenstein algebra of embedding dimension $n, $  it  will be useful to get information on the  homogeneous cubics $F \in P=K[y_1, \dots, y_n] $ such that  $HF_{A_F} =\{1,n,n,1\}, $ that is to characterize the homogeneous cubics which are non degenerate.

\begin{remark} \label{delta} We consider  an homogeneous  form $F_3\in  P=K[y_1, \dots,y_n] $ of degree  three.   Let $\ui\in \mathbb N^n$ be a degree three  multi-index and   write $F_3$ in the dual basis $E^*$
$$
F_3= \sum_{|\ui|=3} \alpha_\ui \; \frac{1}{\ui !} y^\ui .  $$
$F_3 $ is non degenerate ($ HF_{A_{F_3}}(1)=n$) if and only if the $K$-vector space generated by all the derivatives of order  two    has dimension $n, $ that is $$< \partial_{\ui} F_3 \ \ : \ |\ui |=2>_K =P_1. $$

This condition can be formulated in terms of the rank of a matrix, say $ \Delta_{F_3}, $ given by  the coefficients of the linear forms $\partial_{\underline i} F_3 ,  |\ui|=2. $
 The matrix $ \Delta_{F_3}  $ has  size $n \times {{n+1}\choose{2}} $ with entries in the $\alpha_{\underline i}'s  \in K. $   We label the rows by $j =1, \dots, n $ and the columns by $\underline i \in \mathbb N^n$  (a degree two multi-index). We have
\begin{equation} \label{alfa}  (\Delta_{F_3})_{j, \underline i}= \alpha_{\underline i + \delta_j} \end{equation}
where $\delta_j  $ is the $n$-uple with $0$-entries but $1$ in position $j, $ hence $\underline i + \delta_j  =(i_1, \dots, i_j+1, \dots, i_n).$  In fact we have
$$  \partial_{\underline i } F_3= \sum_{|\underline p|=3} \alpha_{\underline p}\  y^{\underline p - \underline i}= \sum_{j =1}^n  \alpha_{\underline i+ \delta_j}\  y_j .$$

\noindent Hence  $F_3$ is  non degenerate   if and only if $rk ( \Delta_{F_3}) =n $.
\end{remark}

 \bigskip
\begin{theorem} \label{n}
Let $A=R/I$ be an Artinian  Gorenstein local ring  with Hilbert function
$\{1,n,n,1\}$.
Then $A $ is canonically graded.
\end{theorem}
\begin{proof}
Let $F=F_0+F_1+F_2+F_3$ be a  polynomial of $P=K[y_1,\dots, y_n] $  of  degree three
such that $I=Ann_R( F)$ ($F_i$ denotes the homogeneous components of  degree $i$).
Since $HF_{A} $ is symmetric, then  $ G= gr_{\m}(A) $ is Gorenstein, in particular  $ G=Q(0) \simeq  R/Ann_R(F_3)  = A_{F_3} $ and $ rk (\Delta_{F_3})= n $ being $F_3$ non degenerate for the Hilbert  function $\{1,n,n,1\}$.
By the admissibility of $F_3$ we deduce that  $ P_{\le 1} \subseteq \langle F_2+F_3 \rangle_R.$ Hence we  may assume $F =F_3 +F_2, $ that is
$$I^\bot=\langle F\rangle_R=\langle F_2+F_3 \rangle_R .$$

So we have to prove that,  however we fix $F_2, $ there exists
   an automorphism $\varphi$ of $R/\mathcal M^4$ which induces
   $$ A_{F_3} \simeq   A_{F_3 + F_2}.$$
Let $\varphi$ be an automorphism of $R/\mathcal M^4$ with the identity as Jacobian defined as follows
$$
\varphi(x_j)=x_j+ \sum_{|\ui|=2} a_{\ui}^j x^{\ui}
$$
for $j=1, \dots, n.$ We prove that there exists $\ua = (a_{\ui}^1, |\ui|=2  ; \cdots; a_{\ui}^n, |\ui|=2 )\in \mathbb N^{n \binom{n+1}{2}}$,  the row vector of the coefficients defining  $\varphi$, such that
\begin{equation}
 \label{start}  [F_3]_{E^*} M(\varphi) = [F]_{E^*}.
 \end{equation}
The  matrix associated to $\varphi, $ say $ M(\varphi),  $  is an  element of $Gl_r(K)$, $r=\binom{n+3}{4}$, with respect to the basis $E$
of $R/\mathcal M^4$ ordered by the deg-lexicographic order,   hence
$$ M(\varphi)=\left(
    \begin{array}{l|l|l|l}
1 & 0& 0 & 0  \\      \hline
0 & I_n& 0 & 0  \\   \hline
0 & D& I_{\binom{n+1}{2}}& 0  \\      \hline
0 & 0& B & I_{\binom{n+2}{3}}
\end{array}
  \right)
$$

\noindent
where for  all $t\ge 1$,
$I_t$ denotes  the $t\times t$ identity matrix.
The first block column corresponds to the image $\varphi(1)=1$;
the second block column corresponds
to the image of $\varphi(x_i)$, $i=1,\dots,n$;
the third block column corresponds to the image of
$\varphi(x^\ui)$ such that $|\ui|=2$;
and finally the last block column corresponds to the image of
$\varphi(x^\ui)$ such that $|\ui|=3$, i.e. the identity matrix.

Hence
$D$ is the $\binom{n+1}{2}\times n$ matrix defined by the coefficients of the degree two monomials of $\varphi(x_i)$,
$i=1,\dots,n $ and $ B$ is a $\binom{n+2}{3}\times \binom{n+1}{2}$ matrix defined by the coefficients of the degree three monomials appearing in
$\varphi(x^\ui), $  $|\ui|=2$.
It is clear  that $M(\varphi) $  is determined by $D$
and  the entries of $B$ are linear forms in the variables
$a_{\ui}^j$, with  $|\ui|=2$, $j=1,\cdots,n$.
Let
$$
F_2= \sum_{|\ui|=2} \beta_\ui \; \frac{1}{\ui !} y^\ui \ \text{ and } F_3=  \sum_{|\ui|=3} \alpha_\ui \; \frac{1}{\ui !} y^\ui
$$
hence   (\ref{start}) is equivalent  to the following equality
$$ [\alpha_\ui ]  B=[\beta_\ui  ].  $$
Then we get a system of ${{n+1}\choose {2}} $ equations which    are bi-homogeneous polynomials in  the $\{\alpha_\ui \} $ and   $\ua \in \mathbb N^{\ n \binom{n+1}{2}} $ of bi-degree $ (1,1). $
Then there exists a matrix $ M_{F_3} $ of size  $  \binom{n+1}{2} \times {\ n \binom{n+1}{2}} $ and  entries   in  the $\{\alpha_\ui \}'s $ such that
$$   [\alpha_\ui ]  B = M_{F_3}  \ ^t \ua $$
where $^ t \ua $ denotes the transpose of the row-vector $\ua.$
We have to prove that the following linear system in  $  \binom{n+1}{2} $ equations and the ${\ n \binom{n+1}{2}} $ indeterminates $ \ua = (a_{\ui}^1,  \cdots; a_{\ui}^n  ) $
$$ M_{F_3} \  ^t \ua =\  ^t [\beta_\ui  ] $$ is compatible.
The result   follows if we  show that   $rk ( M_{F_3}) $ is maximal, i.e. $rk ( M_{F_3})=\binom{n+1}{2}.$

\medskip
\noindent {\bf{Claim.}}
The matrix $M_{F_3}$ has the following upper-diagonal structure
$$ M_{F_3}=\left(
    \begin{array}{l|l|l|l|l}
M^1_{F_3}&  * &  \cdots & * & * \\      \hline
0 & M^2_{F_3}&   \cdots & * & * \\   \hline
\vdots & \vdots&  \vdots & \vdots & \vdots  \\      \hline
0 & 0& 0 &  M^{n-1}_{F_3}&*  \\ \hline
0 & 0& 0 & 0&  M^n_{F_3}
\end{array}
  \right)
$$
where $ M^l_{F_3}$ is a $(n-l+1)\times \binom{n+1}{2}$ matrix, $l =1,\cdots,n$, such that:
\begin{enumerate}
\item[(i)]
$1$-th row of  $M^1_{F_3}$ $=$  $2$  times the 1-th  row of $\Delta_{F_3}$,

$t$-th row of  $M^1_{F_3}$ $=$  $t$-th  row of $\Delta_{F_3}$, $t=2,\cdots,n$.
\item[(ii)]
$1$-th row of $M^l_{F_3}$ $=$ $2 $ times the  $l$-th   row of  $\Delta_{F_3}$,  for $l=2,\cdots,n$,

$t$-th row of $M^l_{F_3}$ $=$ $(l+t-1)$-th  row of  $\Delta_{F_3}$,  for $t=2,\cdots, n-l+1$, $l=2,\cdots,n$,

\end{enumerate}
where  $\Delta_{F_3}$ is the matrix defined in Remark \ref{delta} of the coefficients of the second derivatives of $F_3.$

\medskip
\noindent {\bf Proof of the  Claim. }
Let us recall that the entries of the columns of $B$ are the coefficients of the degree three monomials of $\varphi(x^\ui)$, $|\ui|=2$.
Hence  the entries of the $a_\ui^l-$th  column  of $M_{F_3}$ are the coefficients of the terms of degree three  in the support of $F_3$ which appear in $\varphi(x^\ui)$ with  coefficient $a_\ui^l$. Given integers $1\le l\le j\le n$ let us compute $\varphi(x_l x_j)$.
If $l\neq j$ then
$$
\varphi(x_l x_j)=x_l x_j + \sum_{|\ui|=2} a_\ui^j x^\ui x_l + \sum_{|\ui|=2} a_\ui^l x^\ui x_j +  \text{ terms of degree $4$}.
$$
Since $x^\ui x_j =x^{\ui +\delta_j}$ and $j> l$ we get
$$
(M_{F_3})_{\delta_l+\delta_j, a_\ui^l}= \alpha_{\ui +\delta_j}
$$
and
$$
(M_{F_3})_{\delta_l+\delta_j, a_\ui^j}= \alpha_{\ui +\delta_l}.
$$
If $j = l$ then
$$
\varphi(x_l^2)=x_l^2 + 2 \sum_{|\ui|=2} a_\ui^l x^\ui x_l + \text{ terms of degree $4$}.
$$
so
$$
(M_{F_3})_{2 \delta_l, a_\ui^l}= 2 \alpha_{\ui +\delta_l}.
$$

\noindent Hence  the row $(\delta_l +\delta_j)$-th of $M_{F_3}$, $j> l$, can be split in two non-zero subsets of entries.
The first subset, with respect the lex ordering, corresponds  to the columns  $a_\ui^l$, $|\ui|=2$, with entries $ \alpha_{\ui +\delta_j}$,
the second subset of entries corresponds to
the colons $a_\ui^j$ with entries $ \alpha_{\ui +\delta_l}$.
From these facts  we get the upper-diagonal block structure of $M_{F_3}.$  In particular,  if we fix $l=1, \dots, n, $ the matrices $M^l_{F_3},  $  $l=1, \dots, n $ of the claim are  determined by  the columns  $a_\ui^l$ and the rows
$ \delta_l+\delta_j $ with $ l \le j \le n$ ($(n-j+1)$-rows) and
$$(M^l_{F_3})_{\delta_l+\delta_j, a_\ui^l}= \alpha_{\ui +\delta_j} \ \ \   \text {if\ \ } l>j \ \ \ ; \ \ \   (M^l_{F_3})_{\delta_l+\delta_j, a_\ui^l}=2 \alpha_{\ui +\delta_j} \ \ \   \text {if\ \ } l=j $$
By Remark \ref{delta} (\ref{alfa}) we get
$$(M^l_{F_3})_{\delta_l+\delta_j, a_\ui^l}=(\Delta_{F_3})_{j,\ui} \ \ \   \text {if\ \ } l>j \ \ \ ; \ \ \   (M^l_{F_3})_{\delta_l+\delta_j, a_\ui^l}=2 (\Delta_{F_3})_{j,\ui}  \ \ \   \text {if\ \ } l=j  $$
as claimed.
\noindent

 \medskip

\noindent Now we prove that $$\ \ \ \ \ rk(M_{F_3} ) =     \binom{n+1}{2}. $$

\noindent Since $F_3 $ is non degenerate, by Remark \ref{delta} we have $ rank(\Delta_{F_3}) =n. $ Now $M^l_{F_3} $ for $l=1, \dots, n $ is a matrix of size $n-l+1 \times  \binom{n+1}{2}  $ obtained by $\Delta_{F_3}$  by deleting the first $l $ rows. Hence $rk (M^l_{F_3})= n-l+1 $ and the result follows.
\end{proof}

\vskip 2mm
\noindent From the previous result, we easily get the following consequences. 
\medskip
\begin{corollary} \label{iso}
There exists an isomorphism between the Artinian  Gorenstein local rings   $(A,\m) $ and $(B, \n) $ with Hilbert function
$ \{1,n,n,1\}$ if and only if $gr_{\m}(A) \simeq gr_{\n}(B) $ as graded $K$-algebras.
\end{corollary}

\medskip
\begin{corollary} \label{V(F)}
The classification of   Artinian  Gorenstein local rings
 with Hilbert function  $HF_A=\{1,n,n,1\}$ is equivalent
to the projective classification of the  hypersurfaces
$ V(F)\subset \mathbb P^{n-1}_K$ where  $F$ is  a degree three
non degenerate form in $n$ variables.
\end{corollary}

\bigskip

\bigskip
The classification of the  Artinian Gorenstein  local rings
  with Hilbert function  $HF_A=\{1,n,n,1\}$ for $1 \le n\le 3$ had been studied by Casnati and Notari in \cite{CN09} and by Cartwright et al.  in \cite{CEVV08}. By using Corollary  \ref{V(F)}, the problem is reduced to the homogeneous case which is well known  for $1 \le n\le 3. $ Hence we recover their results  and we describe the geometric models of the varieties defined by them.

\medskip
If $n=1, $ then it is clear that $A\cong K[[x]]/(x^4)$, so there is a only one analytic model.

\medskip
\begin{proposition}
Let $A$ be an Artinian   Gorenstein  local ring with Hilbert function
$HF_A=\{1,2,2,1\}$. Then $A$ is isomorphic to one and only one of the following
  quotients of  $R=K[[x_1,x_2,x_3]]$,$$
\begin{array}{|c|c|c|}           \hline
 \text{ Model } A=R/I& \text{Inverse system } F & \text{Geometry of } C=V(F)\subset \mathbb P^{1}_K \\ \hline
(x_1^3,x_2^2) & y_1^2 y_2 & \text{Double point plus a simple point} \\  \hline
(x_1 x_2, x_1^3 - x_2^3) & y_1^3- y_2^3 & \text{Three distinct  points} \\  \hline
 \end{array}
 $$
\end{proposition}
\begin{proof}
Let us assume $n=2$, then $gr_{\m}(A)=K[y_1,y_2]/Ann(F)$ where $F\in K[y_1,y_2]$ is a degree three form on two
variables $y_1,y_2$.
Since $K$ is an algebraic closed field, $F$ can be decomposed as product of three linear forms
$L_1, L_2, L_3$, i.e. $F=L_1 L_2 L_3$.
We set $d= dim_K \langle L_1, L_2, L_3\rangle$, so  we only have to consider three cases.
If $d=1$ then we can assume $F=y_1^3, $ but in this case
 $F$ is  degenerate.

If $d=2, $ then we can assume  $F=y_1^2 y_2$.
It is easy to see that $Ann(\langle y_1^2 y_2 \rangle)=(x_1^3,x_2^2)$.
If $d=3$ then we can assume $F=y_1^3 - y_2^3$. In this case we get
$Ann(\langle y_1^3 - y_2^3 \rangle)=(x_1 x_2, x_1^3 - x_2^3)$.
Since $V(y_1^2 y_2)$ (resp.  $V(y_1^3 - y_2^3) $) is a  degree three subscheme of $\mathbb P^{1}_K$
with two (resp. three) point basis we get that the algebras of the statement are not isomorphic.
\end{proof}

\medskip
We know that any plane elliptic plane cubic curve $C\subset \mathbb P^{2}_K$ is defined, in a suitable
system of coordinates, by a   Legendre's equation
$$
L_\lambda=y_2^2 y_3- y_1(y_1-y_3)(y_1-\lambda y_3)
$$
with $\lambda \neq 0,1$.
Attached to this  equation we can consider the $j$-invariant
$$
j(\lambda)= 2^8 \; \frac{(\lambda^2-\lambda+1)}{\lambda^2 (\lambda -1 )^2}.
$$
Let us denote
$$
\Sigma_\lambda=\{\lambda, \frac{1}{\lambda} , -\lambda, \frac{1}{1- \lambda}, \frac{\lambda}{\lambda-1},
\frac{\lambda-1}{\lambda} \}.
$$
Notice that $j(\Sigma_\lambda)=j(\lambda)$.
It is well known that
two plane elliptic plane cubic  curves $C_i=V(L_{\lambda_i})\subset \mathbb P^{2}_K$, $i=1,2$, are projectively isomorphic
if and only if $j(\lambda_1)=j(\lambda_2)$, \cite{Har77}.

\bigskip
\begin{proposition}
Let $A$ be an Artinian   Gorenstein local ring  with Hilbert function
$HF_A=\{1,3,3,1\}$. Then $A$ is isomorphic to one and only one of the following
  quotients of  $R=K[[x_1,x_2,x_3]]$,
$$
\begin{array}{|c|c|c|}           \hline
 \text{ Model } A=R/I& \text{Inverse system } F & \text{Geometry of } C=V(F)\subset \mathbb P^{2}_K \\ \hline
 (x_1^2,x_2^2,x_3^2) & y_1 y_2 y_3 & \text{Three independent lines} \\  \hline
 (x_1^2,x_1x_3,x_3x_2^2,x_2^3,x_3^2+x_1x_2) & y_2(y_1y_2- y_3^2) & \text{Conic and a tangent line} \\  \hline
 (x_1^2,x_2^2,x_3^2+6x_1x_2) &  x_3(x_1x_2- x_3^2)& \text{Conic and a non-tangent line} \\  \hline
 (x_3^2, x_1x_2,x_1^2+x_2^2-3 x_1x_3) &  y_2^2 y_3- y_1^2( y_1+ y_3)& \text{Irreducible nodal cubic} \\  \hline
 (x_3^2,x_1x_2, x_1x_3,x_2^3,x_1^3+ 3x_2^2x_3) & y_2^2 y_3 - y_1^3 & \text{Irreducible cuspidal cubic} \\  \hline
 (x_2x_3,x_1x_3,x_1x_2,x_2^3-x_3^3,x_1^3-x_3^3) & y_1^3+y_2^3+y_3^3  &  \text{Elliptic Fermat curve} \\ \hline
 I(\lambda)=(x_1x_2,H_1,H_2) & L_\lambda,\; j(\lambda)\neq 0& \text{Elliptic non Fermat curve} \\ \hline
 \end{array}
 $$
 with
 $H_1=\lambda^2 x_1^2 + \lambda (1+\lambda) x_1 x_3 + (\lambda^2-\lambda+1)x_3^2$,
$H_2=\lambda^2 x_2^2 + \lambda x_1 x_3 + (1+\lambda)x_3^2$,
 and $I(\lambda_1)\cong  I(\lambda_2)$ iff $j(\lambda_1)=j(\lambda_2)$
\end{proposition}
\begin{proof}
Let us assume that $F$ is the product of the  linear forms $l_1,l_2,l_3$.
If  $l_1,l_2,l_3$  are $K$-linear independent we get the first case.
On the contrary, if these linear forms are $K$-linear dependent, we get that $F$ is degenerate.

Let us assume that $F$ is the product of a linear form $l$ and an irreducible
quadric  $Q$. According to the  relative position of  $V(l)$ and $V(Q)$
 we get the second and the  third  case.

Let $F$ be  a degree three irreducible form.
If $C=V(F)$ is singular then  we get the cases fourth and fifth.
If $C=V(F)$ is non-singular, then we may assume that $F=L_\lambda$
for $\lambda\in K\setminus \{0,1\}$, i.e. $C$ is an elliptic cubic curve.
If $j(\lambda)=0 $ then $C$ fits in the orbit of Fermat's curve $y_1^3+y_2^3+y_3^3\equiv 0$ and we get the sixth case.
If $j(\lambda)\neq 0 $ then it is easy to see that
$$
J=(x_1x_2, H_1, H_2) \subset Ann_R(L_\lambda).
$$
Since $\lambda^2-\lambda+1\neq 0$ then $J$ is a complete intersection with Hilbert function $\{1,3,3,1\}$.
Hence $(x_1x_2, H_1, H_2)= Ann_R(L_\lambda)$.
\end{proof}

\begin{remark}
As before, the classification of   Artinian Gorenstein  algebras  with Hilbert function
$\{1, 4, 4, 1\}$ can be obtained by using results on    the classification of the degree three hypersurfaces
of $\mathbb P^{3}$, see for instance \cite{BL98}.
\end{remark}

\bigskip

We present now an unexpected application of Theorem \ref{n}  to the classical problem of the rationality of the Poincar\'e  series. We denote by $P^A_K(z) $ the Poincare' series of a local ring $(A, \m, K), $ that is
$$P^A_K(z)= \sum_{j\ge0} Tor_j^A(K,K) z^j.$$

\noindent We recover  a very recent result by
 Henriques and  Sega, see \cite{HS09} Theorem 4.3, in the particular case of rings.  Among other results,  in the quoted paper the authors studied   the  problem of the rationality in the case of   Artinian Gorenstein local rings  $A  $  such that $\m^4=0  $  and $HF_A(1) \ge 3,  $  under the assumption
that there exists
a non-zero  element $a $ in $\m$ such that $(0: a) $ is a principal ideal. A such element is called an {\it{ exact zero divisor}}. Henriques and Sega proved that if $A$ is Gorenstein of this type, then   the existence of an exact zero divisor assures that the Hilbert function is balanced, that is $HF_A =\{ 1, n, n, 1\} $  (see  \cite{HS09} Proposition 4.2.). One can prove that if $a $ is an exact zero divisor in $A, $ then its initial form $a^* $ is an exact zero divisor in $G $ (see  \cite{HS09}).
Notice that, in the graded case,  if $A=A_F$ with $F$ a generic cubic,
then an exact zero divisor always exists and this implies the existence of  a Koszul filtration (see \cite{CRV01}, Theorem 6.3).
\begin{corollary}[ \cite{HS09} Theorem 4.3.]  \label{ines}  Let $A$ be an    Artinian Gorenstein local ring with $\m^4=0 $ and  $HF_A(1) =n \ge 3.  $  If there exists an exact zero divisor in $A$, then $A$ is Koszul and hence $P^A_K(z) $ is rational.
\end{corollary}
\begin{proof} By  the existence of an exact zero divisor, the Hilbert function of $A$ is balanced, that is $HF_A =\{ 1, n, n, 1\}. $   Hence, by Theorem \ref{n}, $A$ is canonically graded. Because  there exists an exact zero divisor in $G, $ by Remark 3.5 in \cite{HS09} and Proposition 2.3 b) in \cite{CRV01}, we conclude that
$G$ has a Koszul filtration.  As a consequence,   $A$ is Koszul since $G$ is Koszul and $ P^A_K(z) = P^G_K(z)$ is rational.
\end{proof}

\bigskip
\section{Gorenstein Artin algebras with socle degree three}

In \cite{CN09}  Casnati and Notari presented a complete classification of the   Artinian Gorenstein local  algebras with Hilbert functions   $ \{1,m,3,1\}, $ $ m \ge 3.$
In this section we study the Artinian Gorenstein  algebras with Hilbert function $ \{1,m,n,1\}. $
By using a result proved by Iarrobino  (see \cite{Iar94}, Proposition 1.9), a necessary condition for which the numerical function $ \{1,m,n,1\}$   is the Hilbert function of an Artinian Gorenstein  local algebra is that $ m \ge n.$
If $m=n $ we have proved  that  every   Artinian Gorenstein  algebra  with Hilbert function $ \{1,n,n,1\} $ is  canonically graded. This is not longer true if $m > n$  because the Hilbert function is not symmetric.
  In this case  the associated graded ring $G$  is not Gorenstein,  but   another  Gorenstein graded algebra will play the same role:  $ Q(0)  $ (see Section 2.)  which is   the unique Gorenstein graded quotient of $G$ with the same socle degree.   
By Iarrobino's work  we deduce that
$$ HF_{Q(0 )} =  \{1,n,n,1\}. $$
If we deal with  different local rings, we will denote by   $Q_A(0) $ the module   corresponding to the local ring $(A,\m).$

\begin{theorem} \label{MAIN}The following facts are equivalent:
\begin{enumerate}
\item[(a)]  $A$ is an Artinian Gorenstein local ring with Hilbert function $\{1,m,n,1\}, $  $m > n$
\item[(b)]  $ A \simeq A_F $  where $F \in  K[y_1, \dots, y_m], $ $ F=F_3 + y_{n+1}^2 + \dots + y_m^2 $ with  $F_3$  a non degenerate form of degree three in $ K[y_1, \dots, y_n]  $
\end{enumerate}
\end{theorem}

\begin{proof}  


Let $F=F_0+F_1+F_2+F_3$ be a   polynomial of $P=K[y_1,\dots, y_m]$ of degree three
such that $I=Ann_R( F)   $ ($F_i$ is an homogeneous form of degree $i$).
We know  that $Q(0) \simeq R/Ann_R(F_3) $  and it has Hilbert function $\{1,n,n,1\}.$ Since $Q(0) $ is a graded algebra of embedding dimension $n<m$,    there exist  $L_{n+1}, \dots, L_{m} $ independent  linear forms contained in  $ Ann_R(F_3),  $ hence we may assume there exist $L_1, \dots, L_n, \dots, L_m $ generators  in $P_1 $ such that  $ F_3 \in S=K[L_1, \dots, L_n].  $

Since the  Hilbert function of $A$  is $\{1,m,n,1\}, $ and hence $dim_K (I^{\perp})_1=m, $
it is easy to see that
$$\langle F \rangle_R  =\langle  F_2+F_3 \rangle_R .$$
Now we can write $F_2= C  +D $ where $C \in K[L_{n+1}, \dots, L_m] $ and $D $ is  a quadratic form  in the monomials $L_i L_j $  with $ 1\le i \le n, \ 1\le j \le m.$ Since the field $K$ is algebraically closed,     we may assume  there exist  $ \lambda_i \in K $ such that $ F_2= \lambda_m L_m^2 + \dots + \lambda_{n+1} L_{n+1}^2 + D'  $  with $D'$  the corresponding replacement of $D. $ \par  \noindent  Since   $HF_A(1)=m, $ we remark  that,  by (\ref{H1}) and (\ref{H2}),  $ \lambda_m, \dots, \lambda_{n+1} $ are different from zero. Summing up these information, we can conclude that there exists $L_1, \dots, L_n, \dots, L_m $ independent linear forms of $P $ such that
$$ A \simeq A_F $$
where   $F=F_3+  L_m^2 + \dots +  L_{n+1}^2 + H $ with   $F_3 $ an homogeneous  form  of degree three   in $ K[L_1, \dots, L_n]  $  and $H$ an homogeneous form of degree two  in the monomials $L_i L_j $  with $ 1\le i \le n, \ 1\le j \le m.$

Since we are  considering  the linear change of coordinates in $P$ sending  $y_i \to L_i,  $    we  should replace $R$ and $P$ via  the corresponding linear automorphism.  For short,  we still denote  by $R=k[[x_1, \dots, x_m]] $ and   $P=k[y_1, \dots, y_m] $ the corresponding images.
Then we have that
$$ A \simeq A_F $$
where   $F=F_3+  y_{n+1}^2 + \dots +  y_{m}^2 + H $ with   $F_3 $ an homogeneous  form  of degree three   in $ K[y_1, \dots, y_n]  $  and $H$ an homogeneous form of degree two  in the monomials $y_i y_j $  with $ 1\le i \le n, \ 1\le j \le m.$
So we have to prove that,  however we fix $H, $ there exists
   an automorphism $\varphi$ of $R/\mathcal M^4$ which induces
   $$ 
   A_{F_3+  y_{n+1}^2 + \dots +  y_{m}^2} \simeq   A_{F_3+  y_{n+1}^2 + \dots +  y_{m}^2 + H}.
   $$
Let $\ui = (i_1, \dots, i_n) $ be a $n$-multi-index of degree two and let $\varphi$ be the special automorphism of $R/\mathcal M^4$ with the identity as Jacobian defined as follows
$$
\varphi(x_j)=x_j+ \sum_{|\ui|=2} a_{\ui}^j x^{\ui}
$$
for $j=1, \dots, m.  $   We prove that there exists $\ua = (a_{\ui}^1 ; \dots; a_{\ui}^m ) \in  \mathbb N^{m  \binom{n+1}{2}}$  the  vector of the coefficients defining  $ \varphi $ such that
\begin{equation} 
\label{start2}  
[F_3+  y_{n+1}^2 + \dots +  y_{m}^2]_{E^*} M(\varphi) = [F_3+  y_{n+1}^2 + \dots +  y_{m}^2 + H]_{E^*}.
\end{equation}
Repeating  the same computation as in Theorem \ref{n}, the  matrix associated to $\varphi, $ say $ M(\varphi),  $  is an  element of $Gl_r(K)$, $r=\binom{n+3}{4}$, with respect to the basis $E$
of $R/\mathcal M^4$ ordered by the deg-lexicographic order,   hence
$$ M(\varphi)=\left(
    \begin{array}{l|l|l|l}
1 & 0& 0 & 0  \\      \hline
0 & I_m& 0 & 0  \\   \hline
0 & D& I_{\binom{m+1}{2}}& 0  \\      \hline
0 & 0& B & I_{\binom{m+2}{3}}
\end{array}
  \right)
$$

\noindent
Precisely
$D$ is the $\binom{m+1}{2}\times m $ matrix defined by the coefficients of the degree two monomials of $\varphi(x_j)$,
$j=1,\dots,m $  and $ B$ is a $\binom{m+2}{3}\times \binom{m+1}{2}$ matrix defined by the coefficients of the degree three monomials appearing in
$\varphi(x^\ui), $  $|\ui|=2$.
It is clear  that $M(\varphi) $  is determined by $D$
and  the entries of $B$ are linear forms in the variables
$a_{\ui}^j$, with  $|\ui|=2$, $j=1,\cdots,m$. Notice that, by the peculiarity of $\varphi, $ both $D$ and $B$ have several zero-rows (precisely  the rows corresponding to the monomials of degree three divided by $x_i, \ \ i>n $).
Let
\vskip 2mm
 {$
H = \sum_{*} \beta_{i j}  \;   y_i y_j  \ \ $ where { \footnotesize $  i<j, \ 1\le i \le n, 1\le j \le m $ }    
and  
$ F_3=   \sum_{|\ui|=3} \alpha_\ui \; \frac{1}{\ui !} y_1^{i_1} \cdots y_n^{i_n}
$ }

\vskip 2mm \noindent  hence   (\ref{start2}) can be reduced to the following equality
$$ [\alpha_\ui ]  B' =[\beta_{i j}  ]   $$ where $B' $ is the submatrix of $B$ of size $ \binom{n +2}{3} \times [\binom{n+1}{2}+n (m-n) ]$  obtained considering  the  rows and columns corresponding to  the degree three monomials  in $x_1, \dots, x_n$ appearing in $\varphi(x_i x_j) $ with ${i<j, 1\le i \le n, 1\le j \le m}. $
Then we get a system of $\binom{n+1}{2}+n (m-n)  $ equations which    are bi-homogeneous polynomials in  the $\{\alpha_\ui \} $ and   $\ua = (a_{\ui}^1,  \cdots; a_{\ui}^m : \ |\ui|=2 )\in \mathbb N^{\ m \binom{n+1}{2}} $ of bi-degree $ (1,1). $  Then there exists a matrix $ M_{F} $ of size   $[ \binom{n+1}{2}+n (m-n) ]\times {\ m \binom{n+1}{2}} $ and  entries   in  the $\{\alpha_\ui \}'s $ such that
$$   [\alpha_\ui ]  B'  = M_{F}  \ ^t \ua $$

We have to prove that the following linear system in the $  \binom{n+1}{2}+n (m-n) $ equations and ${\ m \binom{n+1}{2}} $ indeterminates $ \ua = (a_{\ui}^1,  \cdots; a_{\ui}^n  ) $:
$$ M_{F } \  ^t \ua =\  ^t [\beta_\ui  ] $$ is compatible. The result follows if we show that   $rk ( M_{F }) $ is maximal, i.e.  $ rk ( M_{F })=\binom{n+1}{2}+n (m-n).$
We prove that the matrix $M_{F }$ has the following upper-diagonal structure
$$ M_{F}=\left(
    \begin{array}{l|l|l|l|l}
M_{F_3}&  * &  \cdots & * & * \\      \hline
0 & \Delta_{F_3}&   \cdots & * & * \\   \hline
\vdots & \vdots&  \vdots & \vdots & \vdots  \\      \hline
0 & 0& 0 & \Delta_{F_3}&*  \\ \hline
0 & 0& 0 & 0&  \Delta_{F_3}
\end{array}
  \right)
\hskip -5mm
\begin{array}{l}
\\
\left.
\begin{array}{l}
\phantom{\vdots}\\
 \phantom{\vdots} \\
 \\
\\
\end{array}
\right\}\,  m-n \text{\, times}
\\
\end{array}
$$ where
 \begin{enumerate}
\item[(i)]
$ M_{F_3}$ is a $  \binom{n+1}{2} \times n \binom{n+1}{2}$ matrix defined in Claim  of Theorem \ref{n}.

\item[(ii)]
$\Delta_{F_3}$ ($m-n $ times) is the $ n \times \binom{n+1}{2} $ matrix defined in Remark \ref{delta} of the coefficients of the second derivatives of $F_3.$

\end{enumerate}

\medskip
Following  the definition  of $B', $ we compute $\varphi(x_i x_j)$, $i<j, 1\le i \le n, 1\le j \le m $.
Hence  the entries of the $a_\ui^l-$th  column  of $M_{F}$ ($l=1, \dots, m$) are the coefficients of the terms of degree three  in the support of $F_3$ which appear in $\varphi(x_i x_j)$  with  coefficient $a_\ui^l$.

If  $1\le i \le j\le n $ we are in the same setting of  Claim   of Theorem \ref{n} and we get $ M_{F_3}$ corresponding to the   $a_\ui^l-$th  columns with $l=1, \dots, n.$

We compute now  $ \varphi(x_i x_j) $    with  $i =1, \dots, n $ and $ j=n+1, \dots, m,  $ then
$$
\varphi(x_i x_j)=x_i x_j + \sum_{|\ui|=2} a_\ui^j x^\ui x_i+ \sum_{|\ui|=2} a_\ui^l x^\ui x_j +  \text{ terms of degree $4$}.
$$
Since $x^\ui x_j  $ does not appear in the support of $F_3$ since $j>n, $ then  for every $j=n+1, \dots, m $ we get
$$
(M_{F})_{\delta_i+\delta_j, a_\ui^j}= \alpha_{\ui +\delta_i}
$$
 Hence
 $$
(M_{F})_{\delta_i+\delta_j, a_\ui^j}= (\Delta_{F_3})_{i,\ui}
$$
and $(i) $ and $ (ii)$  are proved.
\noindent

\bigskip
\noindent Since $F_3 $ is non degenerate for the Hilbert function $\{1, n, n,  1\}, $ then by Remark \ref{delta} $rk(\Delta_{F_3}) =n $ and, from Theorem \ref{n}'s proof,    $   rk(M_{F_3} ) =     \binom{n+1}{2}.  $ 
It follows 
$$rk(M_F) =    \binom{n+1}{2}+n (m-n), $$
as required.
\end{proof}
\medskip

We will extend Corollary \ref{iso} to this more general situation.

\begin{corollary}
 \label{isomn}
There exists an isomorphism between the Artinian Gorenstein  local rings $(A,\m) $ and $(B, \n) $ with Hilbert function
$ \{1,m,n,1\}, $ $ m \ge n, $ if and only if $Q_A(0) \simeq Q_B(0) $ as $K$-algebras.

\end{corollary}
\begin{proof} If $m=n, $ then $Q(0) $ coincide with the associated graded ring and the result follows by Corollary \ref{iso}. Assume $m>n,  $  then  result follows from \thmref{MAIN}   which says that the isomorphism classes of $A_F$ only depend only on $F_3$ and hence on  the  isomorphism classes of $Q(0). $
\end{proof}
\vskip 3mm

\medskip
Next result  extends Corollary \ref{V(F)}.

\begin{corollary}
\label{isomnp}
The classification of    Artinian Gorenstein algebras $A$
 with Hilbert function  $HF_A=\{1,m,n,1\}, $ $m \ge n, $  is equivalent
to the projective classification of the cubic  hypersurfaces
$ V(F)\subset \mathbb P^{n-1}_K$ where  $F$ is  a degree three
non degenerate form in $n$ variables.
\end{corollary}

By taking advantage of  the projective classification of the cubic hypersurfaces in $\mathbb P^{r},  $ with $r \le 3,$ the above result gives  a    complete classification of the   Artinian Gorenstein local  algebras with Hilbert functions   $ \{1,m,n,1\}, $ $ n \le 4.$
\medskip



\providecommand{\bysame}{\leavevmode\hbox to3em{\hrulefill}\thinspace}

\bigskip
\bigskip
\noindent
Juan Elias\\
Departament d'\`Algebra i Geometria\\
Universitat de Barcelona\\
Gran Via 585, 08007 Barcelona, Spain\\
e-mail: {\tt elias@ub.edu}

\bigskip
\noindent
Maria Evelina Rossi\\
Dipartimento di Matematica\\
Universit{\`a} di Genova\\
Via Dodecaneso 35, 16146 Genova, Italy\\
e-mail: {\tt rossim@dima.unige.it}

\end{document}